\documentclass{amsart}
\usepackage{graphicx}
\usepackage{amscd}

\theoremstyle{plain}

\newtheorem{definition}{Definition}

\newtheorem{lemma}{Lemma}

\newtheorem{proposition}{Proposition}

\numberwithin{equation}{section}

\input{tcilatex}

\begin{document}
\title[Planar Harmonic Polynomials]{Planar Harmonic Polynomials of Type B}
\author{Charles F. Dunkl}
\address{Department of Mathematics, Kerchof Hall\\
University of Virginia, Charlottesville, VA 22903}
\email{cfd5z@virginia.edu}
\urladdr{http://www.math.virginia.edu/\symbol{126}cfd5z/}
\date{June 7, 1999}

\begin{abstract}
The hyperoctahedral group acting on $\mathbf{R}^{N}$ is the Weyl group of
type \textit{B }and is associated with a two-parameter family of
differential-difference operators $\{T_{i}:1\leq i\leq N\}$. These operators
are analogous to partial derivative operators. This paper finds all the
polynomials $h$ on $\mathbf{R}^{N}$ which are harmonic, $\Delta _{B}h=0$ and
annihilated by $T_{i}$ for $i>2$, where the Laplacian $\ \Delta
_{B}=\sum_{i=1}^{N}T_{i}^{2}$. They are given explicitly in terms of a novel
basis of polynomials, defined by generating functions. The harmonic
polynomials can be used to find wave functions for the quantum many-body
spin Calogero model.
\end{abstract}

\maketitle

\section{Introduction}

For each finite reflection group there are families of invariant inner
products on the space of polynomials, defined by an algebraic expression,
and by integration with respect to invariant weight functions on the sphere
or on all of Euclidean space. These inner products essentially coincide on
the polynomials harmonic with respect to the associated Laplacian. In this
paper we study certain specific explicit harmonic polynomials associated
with the hyperoctahedral group on $\mathbf{R}^{N}$ . For the reflection
groups on $\mathbf{R}^{2}$ all the harmonic polynomials are known as
expressions in Jacobi polynomials. Also orthogonal bases whose elements are
of generalized Hermite (or Laguerre) type, for the Gaussian weight function
have been determined by means of a construction using nonsymmetric Jack
polynomials. However an orthogonal basis for the weight functions on the
sphere (and the ball or the simplex) has not yet been explicitly found. Here
we consider the analogue of ordinary harmonic polynomials in two variables,
but harmonic for the $N$-variable Laplacian $\Delta _{B}$; that is,
polynomials annihilated by $T_{i}$ for $i>2,$where $\{T_{i}:1\leq i\leq N\}$
is the set of differential-difference operators of type $B$ and $\Delta
_{B}=\sum_{i=1}^{N}T_{i}^{2}.$ In previous work the author introduced a
family of polynomials (the ``$p$-basis'') for which it is easy to write down
polynomials annihilated by any desired subset of $\{T_{i}:1\leq i\leq N\}.$
In this study it is important to select a set of polynomials for which the
harmonic polynomials have ``nice'' coefficients. For example, coefficients
of hypergeometric type (Pochhammer symbols) are ``nice''. We introduce a set
of polynomials which are in the $\mathbf{Q}$-span of the $p$-basis
(coefficients independent of the parameters) and which allow nice
expressions. The definition is by means of generating functions.

For the harmonic polynomials we will find the values at a special point, $%
(1,1,\allowbreak \ldots ,1),$ the leading coefficients, and the $L^{2}$%
-norms; the first two are in closed form using $_{2}F_{1}$ and $_{3}F_{2}$
summations respectively, the last is a balanced $_{4}F_{3}$ sum. Finally
there is a discussion of the important applications of the polynomials,
especially as wave functions for spin Calogero quantum-many-body models.

\section{Overview of Results}

The finite group of orthogonal transformations which is generated by
sign-changes and permutation of coordinates on $\mathbf{R}^{N}$ is called
the Weyl group of type $B$, and will be denoted by $W_{N}$. There is a
family of measures associated with $W_{N}$ in a natural way: for positive
parameters $k,k_{1}$ let 
\begin{equation}
d\mu (x;k,k_{1})=\prod_{i=1}^{N}|x_{i}|^{2k_{1}}\prod_{1\leq i<j\leq
N}|x_{i}^{2}-x_{j}^{2}|^{2k}\exp (-\frac{|x|^{2}}{2})dx  \label{gauss}
\end{equation}
Analysis for functions related to these measures depends on the
differential-difference operators constructed by the author \cite{D1}. Note
that the polynomial terms in $d\mu $ correspond to linear functions
vanishing on the reflecting hyperplanes of the Coxeter group $W_{N}$. The
reflections in $W_{N}$ consist of $\{\sigma _{i}:1\leq i\leq N\}$ and $%
\{\sigma _{ij},\tau _{ij}:1\leq i<j\leq N\}$ defined by 
\begin{eqnarray*}
x\sigma _{i} &=&(x_{1},\ldots ,\overset{i}{-x_{i}},\ldots ,x_{N}), \\
x\sigma _{ij} &=&(x_{1},\ldots ,\overset{i}{x_{j}},\ldots ,\overset{j}{x_{i}}%
,\ldots ,x_{N}), \\
x\tau _{ij} &=&(x_{1},\ldots ,\overset{i}{-x_{j}},\ldots ,-\overset{j}{x_{i}}%
,\ldots ,x_{N}).
\end{eqnarray*}
For notational convenience $\sigma _{ji}=\sigma _{ij}$ and $\tau _{ji}=\tau
_{ij}$. We use the same symbols to indicate the action on functions, for
example, $\sigma _{ij}f(x):=f(x\sigma _{ij}).$ The differential-difference
(``Dunkl'') operators of type $B$ (associated with $W_{N}$) are 
\begin{equation*}
T_{i}:=\frac{\partial }{\partial x_{i}}+k_{1}\frac{1-\sigma _{i}}{x_{i}}%
+k\sum_{j\neq i}\left\{ \frac{1-\sigma _{ij}}{x_{i}-x_{j}}+\frac{1-\tau _{ij}%
}{x_{i}+x_{j}}\right\} ,1\leq i\leq N.
\end{equation*}
They are homogeneous of degree $-1$ on polynomials and commute, $%
T_{i}T_{j}=T_{j}T_{i}$. The Laplacian operator is $\Delta
_{B}:=\sum_{i=1}^{N}T_{i}^{2}.$

In this paper we will determine all polynomials $h$\ on $\mathbf{R}^{N}$%
which satisfy $\Delta _{B}h=0$ and $T_{i}h=0$ for all $i>2$. In the standard
case $k=0=k_{1}$ this implies that $h$ depends only on $x_{1},x_{2}$; this
explains the name ``planar''.

\begin{definition}
\label{FF01}The set of polynomials $\{\phi _{n,j},\psi _{n,j}:0\leq j\leq
n=0,1,2,\ldots \}$ is defined by 
\begin{eqnarray*}
\sum_{n=0}^{\infty }\sum_{j=0}^{n}\phi _{n,j}(x)s^{j}t^{n} &=&F_{0}+F_{1}, \\
\sum_{n=0}^{\infty }\sum_{j=0}^{n}\psi _{n,j}(x)s^{j}t^{n}
&=&x_{1}(F_{0}+sF_{1})-x_{1}F_{1}
\end{eqnarray*}
\ \ \ \ \ in terms of the generating functions ($x\in \mathbf{R}^{N},$ and
absolute convergence holds for $|s|<4/3$ and $|t|<\min (1/x_{i}^{2}:1\leq
i\leq N$ $)/3$): 
\begin{eqnarray*}
F_{0}(x;s,t) &=&\frac{1-st(x_{1}^{2}+x_{2}^{2})+t^{2}x_{1}^{2}x_{2}^{2}}{%
(1-2stx_{1}^{2}+t^{2}x_{1}^{4})(1-2stx_{2}^{2}+t^{2}x_{2}^{4})}%
\prod_{i=1}^{N}(1-2stx_{i}^{2}+t^{2}x_{i}^{4})^{-k}, \\
F_{1}(x;s,t) &=&\frac{t(x_{1}^{2}-x_{2}^{2})}{%
(1-2stx_{1}^{2}+t^{2}x_{1}^{4})(1-2stx_{2}^{2}+t^{2}x_{2}^{4})}%
\prod_{i=1}^{N}(1-2stx_{i}^{2}+t^{2}x_{i}^{4})^{-k}.
\end{eqnarray*}
\end{definition}

In each case, the first and second terms on the right hand side produce the
basis functions with $n+j=0\func{mod}2$ and $n+j=1\func{mod}2$,
respectively. The polynomials $\phi _{n,j}$and $\psi _{n,j}$ are of degrees $%
2n$ and $2n+1$ in $x,$ respectively, for $0\leq j\leq n.$ Because of the
invariance of $\Delta _{B}$ under the subgroup $W_{2}$ acting on the first
two variables, there is a basis of planar harmonic polynomials in which all
the monomials have the same parities of the exponents of $x_{1},x_{2}$; of
course they are even in each of the remaining variables.

\subsection{The harmonic polynomials\label{HP0}:}

The basis elements are labeled $h_{n,0}$ and $h_{n,1}$, where $%
x_{1}^{n}x_{2}^{\varepsilon }$ is the term in $h$ with the highest power of $%
x_{1}$ and $\varepsilon =0$ or $1.$ The formulas depend on the residues$\ 
\func{mod}4.$%
\begin{eqnarray*}
h_{4n,0} &=&\sum_{j=0}^{n}\frac{(2(N-1)k+k_{1}+\frac{1}{2}+2n)_{j}(\frac{1}{2%
})_{j}}{((N-1)k+k_{1}+\frac{1}{2}+n)_{j}(Nk+n+1)_{j}}\phi _{2n,2j}, \\
h_{4n+2,0} &=&\sum_{j=0}^{n}\frac{(2(N-1)k+k_{1}+\frac{3}{2}+2n)_{j}(\frac{1%
}{2})_{j}}{((N-1)k+k_{1}+\frac{3}{2}+n)_{j}(Nk+n+2)_{j}}\phi _{2n+1,2j}, \\
h_{4n+1,1} &=&x_{1}x_{2}\sum_{j=0}^{n}\frac{(2(N-1)k+k_{1}+\frac{3}{2}%
+2n)_{j}(\frac{1}{2})_{j}}{((N-1)k+k_{1}+\frac{3}{2}+n)_{j}(Nk+n+1)_{j}}\phi
_{2n,2j}, \\
h_{4n+3,1} &=&x_{1}x_{2}\sum_{j=0}^{n}\frac{(2(N-1)k+k_{1}+\frac{5}{2}%
+2n)_{j}(\frac{1}{2})_{j}}{((N-1)k+k_{1}+\frac{5}{2}+n)_{j}(Nk+n+2)_{j}}\phi
_{2n+1,2j}
\end{eqnarray*}
The following are of mixed parity: 
\begin{eqnarray*}
h_{4n+1,0} &=&\sum_{j=0}^{n}\frac{(2(N-1)k+k_{1}+\frac{3}{2}+2n)_{j}(\frac{1%
}{2})_{j}}{((N-1)k+k_{1}+\frac{3}{2}+n)_{j}(Nk+n+1)_{j}}\psi _{2n,2j} \\
&&+\sum_{j=1}^{n}\frac{(2(N-1)k+k_{1}+\frac{3}{2}+2n)_{j-1}(\frac{1}{2})_{j}%
}{((N-1)k+k_{1}+\frac{3}{2}+n)_{j-1}(Nk+n+1)_{j}}\psi _{2n,2j-1}
\end{eqnarray*}
and 
\begin{eqnarray*}
h_{4n+3,0} &=&\sum_{j=0}^{n}\frac{(2(N-1)k+k_{1}+\frac{5}{2}+2n)_{j}(\frac{1%
}{2})_{j}}{((N-1)k+k_{1}+\frac{3}{2}+n)_{j}(Nk+n+2)_{j}}\psi _{2n+1,2j} \\
&&+\sum_{j=1}^{n+1}\frac{(2(N-1)k+k_{1}+\frac{5}{2}+2n)_{j-1}(\frac{1}{2}%
)_{j}}{((N-1)k+k_{1}+\frac{3}{2}+n)_{j}(Nk+n+2)_{j-1}}\psi _{2n+1,2j-1}
\end{eqnarray*}
\bigskip

There are two other basis elements, defined as $h_{4n,1}:=\sigma
_{12}h_{4n+1,0}$ and $h_{4n+2,1}:=\sigma _{12}h_{4n+3,0}$. Furthermore, $%
\sigma _{12}h_{4n,0}=h_{4n,0},{}\sigma _{12}h_{4n+1,1}=h_{4n+1,1},$ and $%
\sigma _{12}h_{4n+2,0}=-h_{4n+2,0},{}\sigma _{12}h_{4n+3,1}=-h_{4n+3,1}$
from the obvious symmetry properties of the generating functions.

\subsection{Action of $T_{i}$\label{HP1}:}

Since $\Delta _{B}$ commutes with each $T_{i}$ the action of $T_{1}$ or $%
T_{2}$ on any of the polynomials $h_{n,\varepsilon }$ produces a polynomial
annihilated by $\Delta _{B}$ and $T_{i},i>2,$ that is, a scalar multiple of
another polynomial of this family.\ Specifically the results are (for $%
n=0,1,2,\ldots )$:$\smallskip $

$T_{1}h_{4n,0}=2((N-1)k+n)h_{4n-1,0},\smallskip $

$T_{1}h_{4n+3,0}=-2((N-2)k+k_{1}+n+\frac{1}{2})h_{4n+2,0},\smallskip $

$T_{1}h_{4n+2,0}=2(Nk+n+1)h_{4n+1,0},\smallskip $

$T_{1}h_{4n+1,0}=2((N-1)k+k_{1}+n+\frac{1}{2})h_{4n,0},\smallskip $

$T_{2}h_{4n+3,0}=2(Nk+n+1)h_{4n+1,1},\smallskip $

$T_{2}h_{4n+1,1}=2((N-2)k+k_{1}+n+\frac{1}{2})h_{4n+1,0},\smallskip $

$T_{2}h_{4n+1,0}=2((N-1)k+n)h_{4n-1,1},\smallskip $

$T_{2}h_{4n+3,1}=-2((N-1)k+k_{1}+n+\frac{3}{2})h_{4n+3,0}.\medskip $

These formulae (together with the claim $T_{i}h_{n,\varepsilon }=0$ for $i>2$%
) already imply that $\Delta _{B}h_{n,\varepsilon }=0$ because $%
T_{1}^{2}h_{4n,0}=ch_{4n-2,0}$ for some constant $c$ and $%
T_{2}^{2}h_{4n,0}=\sigma _{12}T_{1}^{2}\sigma _{12}h_{4n,0}=\sigma
_{12}T_{1}^{2}h_{4n,0}=-ch_{4n-2,0}$. \ All the other polynomials $%
h_{m,\varepsilon }$ can be obtained by applying $T_{1}$ or $T_{2}$ often
enough to some $h_{4n,0}$ with $4n>m.$ In subsequent sections we will derive
the values of $||h_{n,\varepsilon }||_{2}^{2},\allowbreak h_{n,\varepsilon
}(1,1,\ldots ,1),\allowbreak h_{n,\varepsilon }(1,0,\ldots ,0)$ and
coefficients of leading terms. All but the $L^{2}$-norms are products of
linear factors in the parameters, while the norms are expressed as sums of
balanced $_{4}F_{3}$-series.

\section{Symbolic Calculus}

The results described above depend on the basis of polynomials introduced in 
\cite{D4} for type $A$, \cite{D5} for type $B.$ The idea is to replace the
variables in the type-$A$ basis by $x_{1}^{2},x_{2}^{2},\ldots ,x_{N}^{2}$
and then use the expressions for $T_{i}$ in terms of corresponding type-$A$
operators. \ Throughout let $y=(y_{1},y_{2},\ldots
,y_{N})=(x_{1}^{2},x_{2}^{2},\ldots ,x_{N}^{2})$ for $x\in \mathbf{R}^{N}$ .
The type-$A$ Dunkl operator is defined by 
\begin{equation*}
\widehat{T}_{i}=\frac{\partial }{\partial y_{i}}+k\sum_{j\neq i}\frac{1-(ij)%
}{y_{i}-y_{j}},1\leq i\leq N.
\end{equation*}
where $(ij)$ denotes the transposition of $y_{i}$ and $y_{j},$ the effect of 
$\sigma _{ij}$ or $\tau _{ij}$ on the squared variables. The polynomials in $%
x\in \mathbf{R}^{N}$ are spanned by polynomials of the form $x^{\varepsilon
}g(y)$ where $\varepsilon =(\varepsilon _{1},\varepsilon _{2},\ldots
,\varepsilon _{N})$ $\ $and $x^{\varepsilon }=x_{1}^{\varepsilon
_{1}}x_{2}^{\varepsilon _{2}}\ldots $ with each $\varepsilon _{i}=0$ or $1.$

\begin{proposition}
\label{Tf0}: Let $f(x)=x^{\varepsilon }g(y)$ with each $\varepsilon _{i}=0$
or $1$. For $i=1,2,\ldots ,N,$%
\begin{eqnarray*}
T_{i}f(x) &=&2x_{i}x^{\varepsilon }\widehat{T}_{i}g(y),\text{ if }%
\varepsilon _{i}=0, \\
T_{i}f(x) &=&2\frac{x^{\varepsilon }}{x_{i}}\left( (k_{1}-\frac{1}{2})g+%
\widehat{T}_{i}(y_{i}g)-k\sum_{j}\{(ij)g:\varepsilon _{j}=1,{}j\neq
i)\}\right) ,\text{ if }\varepsilon _{i}=1.
\end{eqnarray*}
\end{proposition}

This is Proposition 2.1 in \cite{D5}. The $p$-basis for the symmetric group
action is constructed as follows: for $1\leq i\leq N$ the polynomials $%
p_{n}(y_{i};y)$ are given by the generating function 
\begin{equation*}
\sum_{n=0}^{\infty
}p_{n}(y_{i};y)r^{n}=(1-ry_{i})^{-1}\prod_{j=1}^{N}(1-ry_{j})^{-k};
\end{equation*}
then for $\alpha =(\alpha _{1},\ldots ,\alpha _{N})\in \mathbf{Z}_{+}^{N}$
the collection of compositions, the basis element $p_{\alpha
}=\prod_{i=1}^{N}p_{\alpha _{i}}(y_{i};y).$ The key property is that $%
\widehat{T}_{j}p_{n}(y_{i};y)=0$ for $j\neq i$ . It was shown in \cite{D4}
that 
\begin{gather}
\widehat{T}_{i}p_{\alpha }=(Nk+\alpha _{i})p_{\alpha
_{i}-1}(y_{i};y)\prod_{m\neq i}p_{\alpha _{m}}(y_{m};y)+k\sum_{j\neq i}
\label{TPA} \\
\left\{ \sum_{m=0}^{\alpha _{j}-1}(p_{\alpha _{i}+\alpha
_{j}-1-m}(y_{i};y)p_{m}(y_{j};y)-p_{m}(y_{i};y)p_{\alpha _{i}+\alpha
_{j}-1-m}(y_{j};y))\prod_{n\neq i,j}p_{\alpha _{n}}(y_{n};y)\right\} , 
\notag
\end{gather}
if $\alpha _{i}>0,$ and $\widehat{T}_{i}p_{\alpha }=0$ if $\alpha _{i}=0.$
On the right hand side of this formula the same term may appear twice, but
this does not matter for the intended use. We will set up a linear
isomorphism between the span of the $p$-basis and polynomials in the formal
variables $p_{1},p_{2},\ldots ,p_{N}$ which is generally not multiplicative
but does allow a simple formula for $\widehat{T}_{i}.$

\begin{definition}
The linear isomorphism\ $\Psi $ between the span of the $p$-basis and the
space $\mathcal{P}$ of polynomials in the formal variables $%
p_{1},p_{2},\ldots ,p_{N}$ is given by $\Psi :p_{\alpha }\longmapsto
p_{1}^{\alpha _{1}}p_{2}^{\alpha _{2}}\ldots p_{N}^{\alpha _{N}}$ (and
extended by linearity). Further the linear transformations $\zeta _{i,j}$
and $\eta _{i}$ on $\mathcal{P}$ are defined by $\zeta _{i,j}p_{1}^{\alpha
_{1}}p_{2}^{\alpha _{2}}\ldots p_{N}^{\alpha _{N}}=\allowbreak p_{i}^{\alpha
_{i}+\alpha _{j}}\prod_{m\neq i,j}p_{m}^{\alpha _{m}}$ for i$\neq j,$ and $%
\eta _{i}$ is evaluation at $p_{i}=0,$ for $1\leq i\leq N.$ (That is, $\zeta
_{i,j}$ replaces $p_{j}$\ by $p_{i}$, and $\eta _{i}$ replaces $p_{i}$ by $%
0. $)
\end{definition}

It is clear that $\Psi $ commutes with the $S_{N}$-action. \ We use the
simplified notation $\widehat{T}_{i}$ for the operator $\Psi \widehat{T}%
_{i}\Psi ^{-1}$ on $\mathcal{P}.$

\begin{proposition}
For $1\leq i\leq N,$ the operator $\widehat{T}_{i}$ acts on polynomials in $%
\mathcal{P}$ \ by 
\begin{equation*}
\widehat{T}_{i}=\frac{\partial }{\partial p_{i}}+Nk\frac{1-\eta _{i}}{p_{i}}%
+k\sum_{j\neq i}\frac{\zeta _{i,j}+\zeta _{j,i}-1-(i,j)}{p_{i}-p_{j}}.
\end{equation*}
\end{proposition}

\begin{proof}
It suffices to examine the effect of the formula on monomials $p_{1}^{\alpha
_{1}}p_{2}^{\alpha _{2}}\ldots p_{N}^{\alpha _{N}}$ and for $i=1.$ The first
two terms produce $(\alpha _{1}+Nk)$ \ if $\alpha _{1}>0,$ else $0.$ In the
sum, the (typical) term for $j=2$ is $(p_{1}^{\alpha _{1}+\alpha
_{2}}+p_{2}^{\alpha _{1}+\alpha _{2}}-p_{1}^{\alpha _{1}}p_{2}^{\alpha
_{2}}-p_{1}^{\alpha _{2}}p_{2}^{\alpha _{1}})\prod_{m=3}^{N}p_{m}^{\alpha
_{m}}/(p_{1}-p_{2}).$ A simple calculation shows this is the image under $%
\Psi $ of the corresponding term in equation \ref{TPA}.
\end{proof}

We will use generating functions to determine the effects of $T_{1}$ and $%
T_{2}$ on the polynomials defined in the previous sections (now we are
considering the type-$B$ operators). For this purpose we consider the cases $%
f(p),x_{1}f(p),x_{1}x_{2}f(p)$ where $f(p)$ is a formal series in $%
p=(p_{1},p_{2})$ (the validity of term-by-term action comes from the same
argument used to justify term-by-term differentiation of a power series
inside its disk of convergence). In the following, we use the notation: 
\begin{equation*}
\delta
_{1,2}f(p_{1},p_{2})=(f(p_{1},p_{1})+f(p_{2},p_{2})-f(p_{1},p_{2})-f(p_{2},p_{1}))/(p_{1}-p_{2}).
\end{equation*}

\begin{proposition}
For a power series $f(p_{1},p_{2})$ (absolutely convergent in the region $%
\{(p_{1},p_{2}):|p_{1}|<1,|p_{2}|<1)$ \ the following hold:
\end{proposition}

\begin{enumerate}
\item  $T_{1}f(p)=2x_{1}\left( \dfrac{\partial f(p)}{\partial p_{1}}+Nk%
\dfrac{f(p_{1},p_{2})-f(0,p_{2})}{p_{1}}+k\delta _{1,2}f(p)\right) ,$

\item  $T_{1}(x_{1}f(p))=2\left( ((N-1)k+k_{1}+\frac{1}{2})f(p)+p_{1}\dfrac{%
\partial f(p)}{\partial p_{1}}+k\delta _{1,2}(p_{1}f(p))\right) ,$

\item  $T_{2}(x_{1}x_{2}f(p))$\newline
$=2x_{1}\left( ((N-1)k+k_{1}+\frac{1}{2})f(p)+p_{2}\dfrac{\partial f(p)}{%
\partial p_{2}}-k\delta _{1,2}(p_{2}f(p))-kf(p_{2},p_{1})\right) ,$

\item  $T_{2}(x_{1}f(p))=2x_{1}x_{2}\left( \dfrac{\partial f(p)}{\partial
p_{2}}+Nk\dfrac{f(p_{1},p_{2})-f(p_{1},0)}{p_{2}}-k\delta _{1,2}f(p)\right)
. $
\end{enumerate}

\begin{proof}
Formulas (1) and (4) follow immediately from Proposition \ref{Tf0}. It was
shown in Lemma 2.3 of \cite{D5} that $\widehat{T}_{i}y_{i}=\widehat{T}_{i}%
\widehat{\rho }_{i}-k,$ where $\widehat{\rho }_{i}$ is the conjugate under $%
\Psi $ of multiplication by $p_{i}$ acting on $\mathcal{P}$. Together with
Proposition 2 this proves formulas (2) and (3).
\end{proof}

By the fundamental properties of the $p$-basis, $T_{i}(x_{1}^{\varepsilon
_{1}}x_{2}^{\varepsilon _{2}}f)(p_{1},p_{2})=0$ for all $i>2,$ and $%
\varepsilon _{1},\varepsilon _{2}=0$ or $1.$ This applies to all the
polynomials used in the sequel. The images under $\Psi $ of the generating
functions $F_{0},F_{1}$ defined in section 2 are in fact simple rational
functions in $p=(p_{1},p_{2}).$ Indeed for indeterminates $u_{1},u_{2}$,
from the definition of the $p-$basis it follows that 
\begin{gather}
\Psi \left( (1-u_{1}y_{1})^{-1}(1-u_{2}y_{2})^{-1}\prod_{i=1}^{N}\left\{
(1-u_{1}y_{i})(1-u_{2}y_{i})\right\} ^{-k}\right) =  \label{y2p} \\
\Psi (\sum_{m,n=0}^{\infty
}p_{(m,n)}u_{1}^{m}u_{2}^{n})=\sum_{m,n=0}^{\infty
}p_{1}^{m}p_{2}^{n}u_{1}^{m}u_{2}^{n}=(1-u_{1}p_{1})^{-1}(1-u_{2}p_{2})^{-1}
\notag
\end{gather}

The desired expressions result from changing variables to $%
u_{1}=tz,u_{2}=tz^{-1},$ and also $u_{1}=tz^{-1},{}u_{2}=tz,$ and $s=\frac{1%
}{2}(z+z^{-1});$ then the two terms are combined by addition and subtraction
(symmetric and skew-symmetric under $(1,2),$ respectively). To ensure
convergence some region must be chosen, for example, $|s|<\frac{4}{3}$ and $%
|t|<1/(3\max (|p_{1}|,|p_{2}|)$ (and $s,t\in \mathbf{C}$). This is valid
because \TEXTsymbol{\vert}$r-\frac{1}{r}|\leq 2|s|\leq (r+\frac{1}{r})$
where $r=|z|$ and $z\in \mathbf{C},$ thus $|s|<\frac{4}{3}$ implies $\frac{1%
}{3}<|z|<3.$

The method for computing the effect of $T_{i}$ on the polynomials $\phi
_{n,j}$ and $\psi _{n,j}$ is to apply $T_{i}$ to the generating functions
and express the result by means of combinations of $\frac{\partial }{%
\partial s}$ and $\frac{\partial }{\partial t}$ and multiplication by $s,t.$
The two basic functions are: 
\begin{eqnarray*}
w_{1} &=&(1-ztp_{1})^{-1}(1-z^{-1}tp_{2})^{-1} \\
w_{2} &=&(1-z^{-1}tp_{1})^{-1}(1-ztp_{2})^{-1}.
\end{eqnarray*}

Then let 
\begin{eqnarray*}
f_{0} &=&\frac{1}{2}(w_{1}+w_{2})=\frac{1-st(p_{1}+p_{2})+t^{2}p_{1}p_{2}}{%
(1-2stp_{1}+t^{2}p_{1}^{2})(1-2stp_{2}+t^{2}p_{2}^{2})} \\
f_{1} &=&(z-z^{-1})^{-1}(w_{1}-w_{2})=\frac{t(p_{1}-p_{2})}{%
(1-2stp_{1}+t^{2}p_{1}^{2})(1-2stp_{2}+t^{2}p_{2}^{2})}
\end{eqnarray*}

The same formulas apply to the images under $\Psi ^{-1}.$

\begin{proposition}
The generating functions in Definition \ref{FF01} satisy $\Psi F_{0}=f_{0}$
and $\Psi F_{1}=f_{1}$.
\end{proposition}

\begin{proof}
Apply $\Psi ^{-1}$ to $w_{1}$ and $w_{2}$ using equation \ref{y2p}, then
both $\Psi ^{-1}w_{1}$ and $\Psi ^{-1}w_{2}$ have the common factor $%
\prod_{i=1}^{N}\left( (1-u_{1}y_{i})(1-u_{2}y_{i})\right)
^{-k}=\prod_{i=1}^{N}(1-2sty_{i}+t^{2}y_{i}^{2})^{-k}$. The parts of the
calculation involving $(1-u_{i}y_{1})^{-1}$ and $(1-u_{i}y_{2})^{-1},i=1$ or 
$2$ proceed just like those with $w_{1}$ and $w_{2}.$
\end{proof}

\subsection{Action of $T_{i}$ on the generating functions:}

Now we can use the symbolic calculus on $f_{0}$ and $f_{1}$. \ Write $%
g_{0}=f_{0}+sf_{1}$ and $g_{1}=-f_{1}$ for the generating functions for \{$%
\psi _{n,j}\}.$ \ Then $w_{1}=f_{0}+\frac{1}{2}%
(z-z^{-1})f_{1}=g_{0}+z^{-1}g_{1}$ and $w_{2}=f_{0}-\frac{1}{2}%
(z-z^{-1})f_{1}=g_{0}+zg_{1}.$ First the effect of $\delta _{1,2}$ on
various functions is calculated: $\delta _{1,2}f_{0}=tf_{1},{}\delta
_{1,2}f_{1}=0,{}\allowbreak \delta _{1,2}(p_{2}f_{0})=sf_{1},\allowbreak
\delta _{1,2}(p_{2}f_{1})=f_{1},\allowbreak \delta
_{1,2}(g_{0})=tf_{1},{}\delta _{1,2}(g_{1})=0,\allowbreak \delta
_{1,2}(p_{1}g_{0})=0,{}\delta _{1,2}(p_{1}g_{1})=f_{1}.$ These simple
relations are the reason for the use of this particular set of functions.

The differentiations can be done on $w_{1}$ and $w_{2}$ separately. It is
easy to verify that 
\begin{eqnarray*}
\frac{\partial w_{1}}{\partial p_{1}} &=&zt(w_{1}+\frac{t}{2}\frac{\partial
w_{1}}{\partial t})+\frac{z^{2}t}{2}\frac{\partial w_{1}}{\partial z}, \\
\frac{\partial w_{2}}{\partial p_{1}} &=&\frac{t}{z}(w_{2}+\frac{t}{2}\frac{%
\partial w_{2}}{\partial t})-\frac{t}{2}\frac{\partial w_{2}}{\partial z}.
\end{eqnarray*}

Further $p_{1}\dfrac{\partial w_{1}}{\partial p_{1}}=\frac{1}{2}(\dfrac{%
\partial }{\partial t}-\dfrac{\partial }{\partial z})w_{1},\allowbreak p_{2}%
\dfrac{\partial w_{1}}{\partial p_{1}}=\frac{1}{2}(\dfrac{\partial }{%
\partial t}+\dfrac{\partial }{\partial z})w_{2},\allowbreak (1-\eta
_{1})w_{1}/p_{1}=ztw_{1}$ and $(1-\eta _{1})w_{2}/p_{1}=(t/z)w_{2}.$ In the
expressions for $\dfrac{\partial f_{0}}{\partial p_{1}}$ and $\dfrac{%
\partial f_{1}}{\partial p_{1}}$the following equations are used:

$zw_{1}+z^{-1}w_{2}=2sg_{0}+2g_{1},\;z^{2}\dfrac{\partial w_{1}}{\partial z}-%
\dfrac{\partial w_{2}}{\partial z}=2(s^{2}-1)\dfrac{\partial g_{0}}{\partial
s}-2g_{1},\bigskip $

$(z-z^{-1})^{-1}(zw_{1}-z^{-1}w_{2})=g_{0},\;(z-z^{-1})^{-1}(z^{2}\dfrac{%
\partial w_{1}}{\partial z}+\dfrac{\partial w_{2}}{\partial z})=s\dfrac{%
\partial g_{0}}{\partial s}+\dfrac{\partial g_{1}}{\partial s},\bigskip $

$z(\dfrac{\partial w_{1}}{\partial z}-\dfrac{\partial w_{2}}{\partial z}%
)=-2sg_{1}+(1-s^{2})\dfrac{\partial g_{1}}{\partial s},\bigskip $

$z(z-z^{-1})^{-1}(\dfrac{\partial w_{1}}{\partial z}+\dfrac{\partial w_{2}}{%
\partial z})=\dfrac{\partial g_{0}}{\partial s}+g_{1}+s\dfrac{\partial g_{1}%
}{\partial s}.\bigskip $

These, as well as the following equations can be proven by direct
verification (express $w_{1}$ and $w_{2}$ in terms of $g_{0}$ and $g_{1}$,
or $f_{0}$ and $f_{1};$ of course $\frac{\partial }{\partial z}=\frac{1}{2}%
(1-z^{-2})\frac{\partial }{\partial s}$). The formulas are grouped by type:

\subsubsection{Case $T_{1}f(y):\protect\phi \rightarrow \protect\psi $}

\begin{gather*}
T_{1}f_{0}=2x_{1}t(\left[ (Nk+1)s+\frac{st}{2}\frac{\partial }{\partial t}+%
\frac{s^{2}-1}{2}\frac{\partial }{\partial s}\right] g_{0}+ \\
\left[ (N-1)k+\frac{1}{2}+\frac{t}{2}\frac{\partial }{\partial t}\right]
g_{1})
\end{gather*}
\begin{equation*}
T_{1}f_{1}=2x_{1}t\left( \left[ (Nk+1)+\frac{t}{2}\frac{\partial }{\partial t%
}+\frac{s}{2}\frac{\partial }{\partial s}\right] g_{0}+\frac{1}{2}\frac{%
\partial }{\partial s}g_{1}\right) ;
\end{equation*}

\subsubsection{Case $T_{2}(x_{1}x_{2}f(y)):x_{1}x_{2}\protect\phi
\rightarrow \protect\psi $}

\begin{gather*}
T_{2}(x_{1}x_{2}f_{0})=2x_{1}(\left[ (N-2)k+k_{1}+\frac{1}{2}+\frac{t}{2}%
\frac{\partial }{\partial t}\right] g_{0}+ \\
\left[ s((N-1)k+k_{1}+1+\frac{t}{2}\frac{\partial }{\partial t})+\frac{%
s^{2}-1}{2}\frac{\partial }{\partial s}\right] g_{1})
\end{gather*}
\begin{equation*}
T_{2}(x_{1}x_{2}f_{1})=2x_{1}\left( -\frac{1}{2}\frac{\partial }{\partial s}%
g_{0}-\left[ (N-1)k+k_{1}+1+\frac{t}{2}\frac{\partial }{\partial t}+\frac{s}{%
2}\frac{\partial }{\partial s}\right] g_{1}\right) ;
\end{equation*}

\subsubsection{Case $T_{1}(x_{1}f(y)):\protect\psi \rightarrow \protect\phi $%
}

\begin{gather*}
T_{1}(x_{1}g_{0})=2(\left[ (N-1)k+k_{1}+\frac{1}{2}+\frac{t}{2}\frac{%
\partial }{\partial t}+\frac{s}{2}\frac{\partial }{\partial s}\right] f_{0}+
\\
\left[ s((N-1)k+k_{1}+1+\frac{t}{2}\frac{\partial }{\partial t})+\frac{%
s^{2}-1}{2}\frac{\partial }{\partial s}\right] f_{1})
\end{gather*}
\begin{equation*}
T_{1}(x_{1}g_{1})=2\left( -\frac{1}{2}\frac{\partial }{\partial s}f_{0}-%
\left[ (N-2)k+k_{1}+\frac{1}{2}+\frac{t}{2}\frac{\partial }{\partial t}%
\right] f_{1}\right) ,
\end{equation*}

\subsubsection{Case $T_{2}(x_{1}f(y)):\protect\psi \rightarrow x_{1}x_{2}%
\protect\phi $}

\begin{equation*}
T_{2}(x_{1}g_{0})=2x_{1}x_{2}t\left( -\frac{1}{2}\frac{\partial }{\partial s}%
f_{0}+\left[ (N-1)k+1+\frac{t}{2}\frac{\partial }{\partial t}\right]
f_{1}\right) ,
\end{equation*}
\begin{gather*}
T_{2}(x_{1}g_{1})=2x_{1}x_{2}t(\left[ (Nk+1)+\frac{t}{2}\frac{\partial }{%
\partial t}+\frac{s}{2}\frac{\partial }{\partial s}\right] f_{0}- \\
\left[ s(Nk+\frac{3}{2}+\frac{t}{2}\frac{\partial }{\partial t})+\frac{1}{2}%
(s^{2}-1)\frac{\partial }{\partial s}\right] f_{1}).
\end{gather*}

\subsection{Action of $T_{i}$ on basis polynomials\label{tbp}:}

Perusal of these formulas reveals that terms involving $k$ are almost ``on
the diagonal'', that is, setting $k=0$ does not noticeably simplify the
formulas. In the $\{\phi _{n,j},\psi _{n,j}\}$basis, the $k\neq 0$ case is
not any more complicated than $k=0$. This is the advantage of these
polynomials over the ordinary $x$-basis. In each of the formulas, the result
of expanding the equations in $\{\phi _{n,j}\}$ for $f_{i},$ $\{\psi
_{n,j}\} $ for $g_{i}$, ($i=0$ or $1)$ and matching up coefficients of $%
s^{j}t^{n}$ on both sides leads to the following (grouped by the parity of $%
n+j$):

\subsubsection{$n+j=0\;\mathbf{\func{mod}{}}2$}

\begin{eqnarray*}
T_{1}\phi _{n,j} &=&(2Nk+n+j)\psi _{n-1,j-1}+ \\
&&(2(N-1)k+n)\psi _{n-1,j}-(j+1)\psi _{n-1,j+1}, \\
T_{1}\psi _{n,j} &=&(2(N-1)k+2k_{1}+n+j+1)(\phi _{n,j}+\phi
_{n,j-1})-(j+1)\phi _{n,j+1}, \\
T_{2}(x_{1}x_{2}\phi _{n,j}) &=&(2(N-2)k+2k_{1}+n+1)\psi _{n,j}+ \\
&&(2(N-1)k+2k_{1}+n+j+1)\psi _{n,j-1}-(j+1)\psi _{n,j+1}, \\
T_{2}\psi _{n,j} &=&x_{1}x_{2}\left( (2(N-1)k+n+1)\phi _{n-1,j}-(j+1)\phi
_{n-1,j+1}\right) ;
\end{eqnarray*}

\subsubsection{$n+j=1{}\func{mod}{}2$}

\begin{eqnarray*}
T_{1}\phi _{n,j} &=&(2Nk+n+j+1)\psi _{n-1,j}+(j+1)\psi _{n-1,j+1}, \\
T_{1}\psi _{n,j} &=&-(2(N-2)k+2k_{1}+n+1)\phi _{n,j}-(j+1)\phi _{n,j+1}, \\
T_{2}(x_{1}x_{2}\phi _{n,j}) &=&-(2(N-1)k+2k_{1}+n+j+2)\psi _{n,j}-(j+1)\psi
_{n,j+1}, \\
T_{2}\psi _{n,j} &=&x_{1}x_{2}\left( (2Nk+n+j+1)(\phi _{n-1,j}-\phi
_{n,j-1})+(j+1)\phi _{n-1,j+1}\right) .
\end{eqnarray*}

Notice that each expression has no more than three different polynomials on
the right-hand side. In the next section we use these to determine the
harmonic polynomials.

\section{Properties of the harmonic polynomials}

Here we demonstrate the action of $T_{1},T_{2}$ on the harmonic polynomials,
which suffices to show $\Delta _{B}h_{n,\varepsilon }=0$ for each such
polynomial, as mentioned before. By construction the polynomials defined in
section \ref{HP0} satisfy $T_{i}h_{n,\varepsilon }=0$ for $i>2;$ it suffices
to establish the formulae of section \ref{HP1}. In a sense the proofs depend
on induction. Since the computations of $T_{i}h_{n,\varepsilon }$ are
somewhat repetitive we will not give details on each formula.\ The
calculations are direct; the definitions of $F_{0},F_{1}$ and $%
h_{n,\varepsilon }$ were formulated after computer-algebra-aided
experimentation. In addition we determine the values at $(1,1,\ldots
,1),\allowbreak (1,0,\ldots ,0)$ and the $L^{2}$-norms.

\subsection{The action of $T_{i}$ on $h$:}

There are four main cases, each with two parts. Because $k_{1}$ always
appears in the same way we introduce an (abbreviation) notation: 
\begin{equation*}
k_{2}=(N-1)k+k_{1}+\frac{1}{2}.
\end{equation*}

\subsubsection{Case: $T_{1}:h_{2m,0}\rightarrow h_{2m-1,0}$}

We begin with the proof of $T_{1}h_{4n,0}=2((N-1)k+n)h_{4n-1,0},$ for $n\geq
1.$ We write $T_{1}\sum_{i=0}^{m}c_{i}\phi
_{m,i}=\sum_{i=0}^{m-1}(T_{1}^{\ast }c)_{i}\psi _{m-1,i}$ , then 
\begin{eqnarray*}
(T_{1}^{\ast }c)_{i} &=&(2(N-1)k+m)c_{i}+ic_{i-1},\text{ for }m+i=0\func{mod}%
2, \\
(T_{1}^{\ast }c)_{i} &=&(2Nk+m+i+1)(c_{i+1}+c_{i})-ic_{i-1},\text{for }m+i=1%
\func{mod}2.
\end{eqnarray*}

These follow from the equations in section \ref{tbp}. Now let $%
h_{4n,0}=\sum_{j=0}^{n}a_{j}\phi _{2n,2j}$ with\medskip\ $a_{j}=\dfrac{%
((N-1)k+k_{2}+2n)_{j}(1/2)_{j}}{(k_{2}+n)_{j}(Nk+n+1)_{j}}.$ Use the above
equations with $m=2n,$then $T_{1}h_{4n,0}=\medskip \allowbreak
\sum_{j=0}^{n-1}b_{j}\psi _{2n-1,2j}+\allowbreak \sum_{j=1}^{n}c_{j}\psi
_{2n-1,2j-1}$ with $b_{j}=2((N-1)k+n)a_{j}$ and 
\begin{eqnarray*}
c_{j} &=&2(Nk+n+j)a_{j}-2(j-\frac{1}{2})a_{j-1} \\
&=&2((N-1)k+n)\dfrac{((N-1)k+k_{2}+2n)_{j-1}(\frac{1}{2})_{j}}{%
(k_{2}+n)_{j}(Nk+n+1)_{j-1}},
\end{eqnarray*}
which is the claimed\medskip\ result. Similarly one can show $%
T_{1}h_{4n+2,0}=\allowbreak 2(Nk+n+1)h_{4n+1,0}$.

\subsubsection{Case : $T_{1}:h_{2m+1,0}\rightarrow h_{2m,0}$}

Next we consider the case $T_{1}\sum_{i=0}^{m}c_{i}\psi
_{m,i}=\sum_{i=0}^{m}(T_{1}^{\ast }c)_{i}\phi _{m,i},$ where (see section 
\ref{tbp}) 
\begin{eqnarray*}
(T_{1}^{\ast }c)_{i} &=&(2k_{2}+m+i)c_{i}-ic_{i-1},\text{ for }m+i=0\func{mod%
}2, \\
(T_{1}^{\ast }c)_{i} &=&(2k_{2}+m+i+1)c_{i+1}-(2k_{2}-2k+m)c_{i}-ic_{i-1},%
\text{for }m+i=1\func{mod}2.
\end{eqnarray*}

Write 
\begin{equation*}
h_{4n+3,0}=\sum_{j=0}^{n}b_{j}\psi _{2n+1,2j}+\sum_{j=1}^{n+1}c_{j}\psi
_{2n+1,2j-1},
\end{equation*}
then\smallskip\ 
\begin{equation*}
T_{1}h_{4n+3,0}=\sum_{j=0}^{n}a_{j}\phi _{2n+1,2j}+\sum_{j=1}^{n+1}d_{j}\phi
_{2n+1,2j-1},
\end{equation*}
with $d_{j}=2(k_{2}+n+j)c_{j}-(2j-1)b_{j-1}=0$ and 
\begin{eqnarray*}
a_{j} &=&2(k_{2}+n+j+1)c_{j+1}-(2k_{2}-2k+2n+1)b_{j}-2jc_{j} \\
&=&-2((N-2)k+k_{1}+n+\frac{1}{2})\dfrac{((N-1)k+k_{2}+2n+1)_{j}(\frac{1}{2}%
)_{j}}{(k_{2}+n+1)_{j}(Nk+n+2)_{j}},
\end{eqnarray*}
$\smallskip \allowbreak $the claimed multiple of $h_{4n+2,0}.$ \
Similarly\medskip\ write 
\begin{eqnarray*}
h_{4n+1,0} &=&\sum_{j=0}^{n}b_{j}\psi _{2n,2j}+\sum_{j=1}^{n}c_{j}\psi
_{2n,2j-1}, \\
T_{1}h_{4n+1,0} &=&\sum_{j=0}^{n}a_{j}\phi _{2n,2j}+\sum_{j=1}^{n}d_{j}\phi
_{2n,2j-1},
\end{eqnarray*}
with $a_{j}=2(k_{2}+n+j)b_{j}-2jc_{j}$ which is $2(k_{2}+n)$ times the
corresponding coefficient of $h_{4n,0},$ while\ $d_{j}=\allowbreak
2(k_{2}+n+j)b_{j}-2(k_{2}-k+n)c_{j}-\allowbreak (2j-1)b_{j-1}=0.$

\subsubsection{Case : $T_{2}:h_{2m+1,0}\rightarrow h_{2m-1,1}$}

For $T_{2}\sum_{i=0}^{m}c_{i}\psi _{m,i}=\allowbreak
\sum_{i=0}^{m-1}(T_{2}^{\ast }c)_{i}(x_{1}x_{2}\phi _{m-1,i})$ one has 
\begin{gather*}
(T_{2}^{\ast }c)_{i}=(2(N-1)k+m+1)c_{i}-(2Nk+m+i+2)c_{i+1}+ic_{i-1}, \\
\text{ for }m+i=0\func{mod}2, \\
(T_{2}^{\ast }c)_{i}=(2Nk+m+i+1)c_{i}-ic_{i-1},\text{ for }m+i=1\func{mod}2.
\end{gather*}

\subsubsection{Case : $T_{2}:h_{2m+1,1}\rightarrow h_{2m+1,0}$}

For $T_{2}\sum_{i=0}^{m}c_{i}x_{1}x_{2}\phi
_{m,i}=\sum_{i=0}^{m}(T_{2}^{\ast }c)_{i}\psi _{m,i}$ one has 
\begin{eqnarray*}
(T_{2}^{\ast }c)_{i} &=&(2k_{2}-2k+m)c_{i}-ic_{i-1},\text{ for }m+i=0\func{%
mod}2, \\
(T_{2}^{\ast }c)_{i} &=&(2k_{2}+m+i+1)(c_{i+1}-c_{i})-ic_{i-1},\text{ for }%
m+i=1\func{mod}2.
\end{eqnarray*}

\subsection{Values at $(1,1,\ldots ,1)$}

Substituting $x=1^{N}$ $=\allowbreak (1,1,\ldots ,1)\in \mathbf{R}^{N}$  in $%
F_{0}$ and $F_{1}$ produces 
\begin{eqnarray*}
\sum_{n=0}^{\infty }\sum_{j=0}^{n}\phi _{n,j}(1^{N})s^{j}t^{n}
&=&(1-2st+t^{2})^{-(Nk+1)} \\
&=&\sum_{m,i}t^{m}s^{m-2i}2^{m-2i}(-1)^{i}\frac{(Nk+1)_{m-i}}{i!(m-2i)!},
\end{eqnarray*}
since only terms with $n+j=0\func{mod}2$ can have non-zero values; note that
the second equality is familiar as the generating function for Gegenbauer
polynomials. To derive this expansion, write 
\begin{eqnarray*}
(1-2st+t^{2})^{-(Nk+1)} &=&(1+t^{2})^{-(Nk+1)}(1-\frac{2st}{1+t^{2}}%
)^{-(Nk+1)} \\
&=&\sum_{j=0}^{\infty }\frac{(Nk+1)_{j}}{j!}(2st)^{j}(1+t^{2})^{-(Nk+1+j)} \\
&=&\sum_{j=0}^{\infty }\sum_{i=0}^{\infty }\dfrac{(Nk+1)_{j+i}}{j!i!}%
(2s)^{j}(-1)^{i}t^{2i+j},
\end{eqnarray*}
now let $j=m-2i.$

\begin{proposition}
For $n=0,1,2,\ldots $ the following hold: 
\begin{eqnarray*}
h_{4n,0}(1^{N}) &=&\frac{(Nk+1)_{n}((N-1)k+1)_{n}}{n!(k_{2}+n)_{n}}, \\
h_{4n+1,0}(1^{N}) &=&\frac{(Nk+1)_{n}((N-1)k+1)_{n}}{n!(k_{2}+n+1)_{n}}, \\
h_{4n+2,0}(1^{N}) &=&h_{4n+3,1}(1^{N})=0, \\
h_{4n+1,1}(1^{N}) &=&\frac{(Nk+1)_{n}((N-1)k+1)_{n}}{n!(k_{2}+n+1)_{n}}, \\
h_{4n+3,0}(1^{N}) &=&\frac{(Nk+1)_{n+1}((N-1)k+1)_{n}}{n!(k_{2}+n+1)_{n+1}}.
\end{eqnarray*}

\begin{proof}
The non-zero cases are all $_{2}F_{1}$ summations. For the first case, $\phi
_{2n,2j}(1^{N})=\dfrac{2^{2j}(Nk+1)_{n+j}(-1)^{n-j}}{(n-j)!(2j)!}=\dfrac{%
(-1)^{n}(-n)_{j}(Nk+1)_{n+j}}{n!j!(\frac{1}{2})_{j}},$ thus 
\begin{eqnarray*}
h_{4n,0}(1^{N}) &=&\frac{(-1)^{n}(Nk+1)_{n}}{n!}\sum_{j=0}^{n}\frac{%
(-n)_{j}(k_{2}+(N-1)k+2n)_{j}}{(k_{2}+n)_{j}j!} \\
&=&\frac{(-1)^{n}(Nk+1)_{n}(-(N-1)k-n)_{n}}{n!(k_{2}+n)_{n}}.
\end{eqnarray*}
This uses the Chu-Vandermonde sum $_{2}F_{1}\left( \QDATOP{-n,b}{c};1\right)
=\dfrac{(c-b)_{n}}{(c)_{n}};$ and $(-b-n)_{n}=(-1)^{n}(b+1)_{n}$ for
arbitrary $b,c$ and $n=0,1,2,...$. The other formulas are proved in the same
way.
\end{proof}
\end{proposition}

\subsection{Leading Coefficients\label{LC}}

Let $cof(f,x_{1}^{m}x_{2}^{n})$ denote the coefficient of the monomial $%
x_{1}^{m}x_{2}^{n}$ in the expansion of the polynomial $f$ in terms of $%
x_{1},x_{2},\ldots ,x_{N}.$ We will determine the values of $%
cof(h_{m,\varepsilon },x_{1}^{m}x_{2}^{\varepsilon })$. For the case $%
\varepsilon =0$ these values agree with evaluation at $(1,0,\ldots ,0).$
Thus evaluate $F_{0}$ and $F_{1}$ at this point to obtain $%
(1-st)(1-2st+t^{2})^{-(k+1)}=a_{0}$ and $t(1-2st+t^{2})^{-(k+1)}=a_{1},$
respectively. The term $a_{1}$ is multiplied by $(-1)$ to obtain
coefficients of $x_{2}^{n}$ in $h_{n,0}.$ By expansion methods similar to
the previous we obtain: 
\begin{eqnarray*}
a_{0} &=&\sum_{m,i}t^{m}s^{m-2i}2^{m-2i-1}(-1)^{i}\frac{(k+1)_{m-1-i}(2k+m)}{%
i!(m-2i)!}, \\
a_{1} &=&\sum_{m,i}t^{m+1}s^{m-2i}2^{m-2i}(-1)^{i}\frac{(k+1)_{m-i}}{%
i!(m-2i)!}.
\end{eqnarray*}

\subsubsection{Case ($4n+\protect\varepsilon ,\protect\varepsilon )$}

The computation for $h_{4n,0}(1,0,...)$ proceeds as follows: $\phi
_{2n,2j}(1,0,...)=\allowbreak \allowbreak (-1)^{n}\dfrac{%
(-n)_{j}(k+1)_{n}(k+n)_{j}}{n!j!(\frac{1}{2})_{j}}$ and so 
\begin{eqnarray*}
h_{4n,0}(1,0,...) &=&\dfrac{(-1)^{n}(k+1)_{n}}{n!}{}_{3}F_{2}\left( \QDATOP{%
-n,k+n,k_{2}+(N-1)k+2n}{k_{2}+n,Nk+n+1};1\right) \\
&=&\dfrac{(-1)^{n}(k+1)_{n}(k_{2}-k)_{n}((N-1)k+1)_{n}}{%
n!(k_{2}+n)_{n}(Nk+n+1)_{n}}.
\end{eqnarray*}
Clearly $%
cof(h_{4n,0},x_{1}^{4n})=cof(h_{4n,0},x_{2}^{4n})=h_{4n,0}(1,0,...). $ The
sum is an application of the Saalsch\"{u}tz formula $_{3}F_{2}\left( \QDATOP{%
-n,a,b}{c,d};1\right) =\dfrac{(c-a)_{n}(d-a)_{n}}{(c)_{n}(d)_{n}},$ provided 
$-n+a+b+1=c+d.$ \medskip

The corresponding formulae for $h_{4n+1,1}$ are obtained by merely
incrementing $k_{1}$ (and also $k_{2})$ by $1$ $.$ Thus 
\begin{eqnarray*}
cof(h_{4n+1,1},x_{1}^{4n+1}x_{2}) &=&cof(h_{4n+1,1},x_{1}x_{2}^{4n+1}) \\
&=&\dfrac{(-1)^{n}(k+1)_{n}(k_{2}+1-k)_{n}((N-1)k+1)_{n}}{%
n!(k_{2}+n+1)_{n}(Nk+n+1)_{n}}.
\end{eqnarray*}

\subsubsection{Case $(4n+2+\protect\varepsilon ,\protect\varepsilon )$}

For $h_{4n+2,0}(1,0,...)$ we begin with $\phi
_{2n+1,2j}(1,0,...)=\allowbreak (-1)^{n}\dfrac{(-n)_{j}(k+1)_{n+j}}{n!j!(%
\frac{1}{2})_{j}}$ and 
\begin{eqnarray*}
h_{4n+2,0}(1,0,...) &=&\dfrac{(-1)^{n}(k+1)_{n}}{n!}{}_{3}F_{2}\left( 
\QDATOP{-n,k+n+1,k_{2}+(N-1)k+2n+1}{k_{2}+n+1,Nk+n+2};1\right) \\
&=&\dfrac{(-1)^{n}(k+1)_{n}(k_{2}-k)_{n}((N-1)k+1)_{n}}{%
n!(k_{2}+n+1)_{n}(Nk+n+2)_{n}}.
\end{eqnarray*}
Further $%
cof(h_{4n+2,0},x_{1}^{4n+2})=-cof(h_{4n+2,0},x_{2}^{4n+2})=h_{4n+2,0}(1,0,...). 
$ Replace $k_{2}$ by $k_{2}+1$ to obtain the value of $%
cof(h_{4n+3,1},x_{1}^{4n+3}x_{2})=-cof(h_{4n+3,1},x_{1}x_{2}^{4n+3}).$

\subsubsection{Case $(4n+1,0)$}

For $h_{4n+1,0}$ note that $\psi _{2n,2j}=x_{1}(\phi _{2n,2j}+\phi
_{2n,2j-1})$ and $\psi _{2n,2j-1}=-x_{1}\phi _{2n,2j-1}.$ To sketch the
argument, let $h_{4n+1,0}=\sum_{j=0}^{n}a_{j}\psi
_{2n,2j}+\sum_{j=1}^{n}b_{j}\psi _{2n,2j-1}$ and let $\phi
_{2n,2j}(1,0,...)=c_{j}$ and $\phi _{2n,2j-1}(1,0,...)=d_{j}.$ Then $%
\sum_{j=0}^{n}a_{j}c_{j}+\varepsilon \sum_{j=1}^{n}(a_{j}-b_{j})d_{j}$
equals $cof(h_{4n+1,0},x_{1}^{4n+1})$ when $\varepsilon =1,$ and \newline
$cof(h_{4n+1,0},x_{1}x_{2}^{4n})$ when $\varepsilon =-1.$ The value of $%
c_{j} $ was found above, and thus the first sum equals\smallskip\ $\dfrac{%
(-1)^{n}(k+1)_{n}(k_{2}-k+1)_{n}((N-1)k+1)_{n}}{n!(k_{2}+n+1)_{n}(Nk+n+1)_{n}%
}.$ For the second sum, 
\begin{equation*}
a_{j}-b_{j}=((N-1)k+n)\dfrac{(k_{2}+(N-1)k+2n+1)_{j-1}(\frac{1}{2})_{j}}{%
(k_{2}+n+1)_{j}(Nk+n+1)_{j}}.
\end{equation*}
The value of $d_{j}$ is calculated similarly to $\phi _{2n+1,2j}$ and thus
the second sum equals 
\begin{gather*}
\frac{((N-1)k+n)(-1)^{n-1}(k+1)_{n}}{(n-1)!(k_{2}+n+1)(Nk+n+1)}{}%
_{3}F_{2}\left( \QDATOP{1-n,k+n+1,k_{2}+(N-1)k+2n+1}{k_{2}+n+2,Nk+n+2}%
;1\right) \\
=\dfrac{(-1)^{n-1}(k+1)_{n}(k_{2}-k+1)_{n-1}((N-1)k+1)_{n}}{%
(n-1)!(k_{2}+n+1)_{n}(Nk+n+1)_{n}}.
\end{gather*}
Combining the two sums and setting $\varepsilon =1$ and -1 respectively we
obtain 
\begin{equation*}
cof(h_{4n+1,0},x_{1}^{4n+1})=\dfrac{%
(-1)^{n}(k+1)_{n}(k_{2}-k)_{n}((N-1)k+1)_{n}}{n!(k_{2}+n+1)_{n}(Nk+n+1)_{n}}
\end{equation*}
and $cof(h_{4n+1,0},x_{1}x_{2}^{4n})=cof(h_{4n+1,0},x_{1}^{4n+1})\dfrac{%
k_{2}-k+2n}{k_{2}-k}.$

\subsubsection{Case\label{c4n3} $(4n+3,0)$}

For $h_{4n+3,0}$ note that $\psi _{2n+1,2j}=-x_{1}\phi _{2n+1,2j}$ and 
\newline
$\psi _{2n+1,2j-1}=x_{1}(\phi _{2n+1,2j-1}+\phi _{2n+1,2j-2})$. As before,
let $h_{4n+3,0}=\allowbreak \sum_{j=0}^{n}a_{j}\psi
_{2n+1,2j}+\sum_{j=1}^{n+1}b_{j}\psi _{2n+1,2j-1}$ (not the same
coefficients as above) and let $\phi _{2n+1,2j-1}(1,0,...)=c_{j}$ and $\phi
_{2n+1,2j}(1,0,...)=d_{j}.$ Then $\allowbreak \varepsilon
\sum_{j=0}^{n}(-a_{j}+b_{j+1})d_{j}+\sum_{j=1}^{n+1}b_{j}c_{j}$ equals $%
cof(h_{4n+3,0},x_{1}^{4n+3})$ when $\varepsilon =1$ and $%
cof(h_{4n+3,0},x_{1}x_{2}^{4n+2})$ when $\varepsilon =-1.$ By a calculation
similar to the previous one the first sum is found to equal 
\begin{equation*}
(k_{2}+n+\frac{1}{2})\dfrac{(-1)^{n+1}(k+1)_{n}(k_{2}-k+1)_{n}((N-1)k+1)_{n}%
}{n!(k_{2}+n+1)_{n+1}(Nk+n+2)_{n}},
\end{equation*}
and the second sum is 
\begin{equation*}
(k+n+\frac{1}{2})\dfrac{(-1)^{n}(k+1)_{n}(k_{2}-k+1)_{n}((N-1)k+1)_{n}}{%
n!(k_{2}+n+1)_{n+1}(Nk+n+2)_{n}}
\end{equation*}
Combining the two sums and setting $\varepsilon =1$ and -1 respectively we
obtain 
\begin{equation*}
cof(h_{4n+3,0},x_{1}^{4n+3})=\dfrac{%
(-1)^{n+1}(k+1)_{n}(k_{2}-k)_{n+1}((N-1)k+1)_{n}}{%
n!(k_{2}+n+1)_{n+1}(Nk+n+2)_{n}}
\end{equation*}
and $cof(h_{4n+3,0},x_{1}x_{2}^{4n+2})=-cof(h_{4n+3,0},x_{1}^{4n+3})\dfrac{%
k_{2}+k+2n+1}{k_{2}-k}.$

\subsection{Norms:}

For arbitrary polynomials three different inner products have been defined.
However for harmonic polynomials there is really only one. Write 
\begin{equation*}
d\mu _{S}(x;k,k_{1})=\prod_{i=1}^{N}|x_{i}|^{2k_{1}}\prod_{1\leq i<j\leq
N}|x_{i}^{2}-x_{j}^{2}|^{2k}d\omega (x)
\end{equation*}
for the measure on the unit sphere $S=\{x\in \mathbf{R}^{N}:|x|=1\}$, where $%
d\omega $ denotes the normalized roation-invariant surface measure. See
formula (\ref{gauss}) for the definition of $d\mu (x;k,k_{1}).$

\begin{proposition}
Suppose $f,g$ are harmonic ($\Delta _{B}f=0=\Delta _{B}g)$ homogeneous
polynomials of degree $m,n$ respectively. Then 
\begin{eqnarray*}
f(T_{1},T_{2},\ldots )g(x)|_{x=0} &=&c\int_{\mathbf{R}^{N}}f(x)g(x)d\mu
(x;k,k_{1})= \\
&=&\delta _{mn}2^{n}(Nk_{2})_{n}{}c_{S}\int_{S}fgd\mu _{S}.
\end{eqnarray*}
\end{proposition}

This was shown in Theorem 3.8 of \cite{D2}; the normalizing constants
satisfy $c\int_{\mathbf{R}^{N}}d\mu =1=c_{S}\int_{S}d\mu _{S},$ (the
well-known Macdonald-Mehta-Selberg integral). Specializing to the harmonic
polynomials $g=h_{n,\varepsilon },$ for which $T_{2}^{2}g=-T_{1}^{2}g$ we
see that if $f(x)=x_{1}^{\varepsilon _{1}}x_{2}^{\varepsilon
_{2}}f_{0}(x_{1}^{2},x_{2}^{2},\ldots ,x_{N}^{2})$ with $\varepsilon
_{1},\varepsilon _{2}=0$ or 1, and $f_{0}$ is homogeneous of degree $m$ then 
\begin{equation*}
f(T)g(x)=T_{1}^{\varepsilon _{1}}T_{2}^{\varepsilon
_{2}}f_{0}(T_{1}^{2},-T_{1}^{2},0,\ldots ,0)g(x)=f_{0}(1,-1,0,\ldots
,0)T_{1}^{\varepsilon _{1}+2m}T_{2}^{\varepsilon _{2}}g(x)
\end{equation*}
By construction the polynomials $h_{n,\varepsilon }$ are pairwise orthogonal
so only $h_{n,\varepsilon }(T)h_{n,\varepsilon }(x)$ need be computed. We
begin with the calculation of $T_{1}^{n}T_{2}^{\varepsilon }h_{n,\varepsilon
}$. The answers are best stated using a notation introduced in \cite{D5}.

\begin{definition}
For m,n$\in \mathbf{Z}_{+}$ and $m\geq n$ let 
\begin{equation*}
\Lambda (m,n)=(Nk+1)_{a}((N-1)k+1)_{b}(k_{2})_{m-a}(k_{2}-k)_{n-b},
\end{equation*}
where $a=\lfloor \frac{m}{2}\rfloor $ and $b=\lfloor \frac{n}{2}\rfloor .$
\end{definition}

This is a special case of the generalized Pochhammer symbol for two-part
partitions. From the formulae in section \ref{HP1} we have 
\begin{eqnarray*}
T_{1}^{4n}h_{4n,0} &=&2^{4n}(-1)^{n}\Lambda (2n,2n), \\
T_{1}^{4n+1}h_{4n+1,0} &=&2^{4n+1}(-1)^{n}\Lambda (2n+1,2n), \\
T_{1}^{4n+2}h_{4n+2,0} &=&2^{4n+2}(-1)^{n}\Lambda (2n+2,2n), \\
T_{1}^{4n+1}T_{2}h_{4n+1,1} &=&2^{4n+2}(-1)^{n}\Lambda (2n+1,2n+1), \\
T_{1}^{4n+3}h_{4n+3,0} &=&2^{4n+3}(-1)^{n+1}\Lambda (2n+2,2n+1), \\
T_{1}^{4n+3}T_{2}h_{4n+3,1} &=&2^{4n+4}(-1)^{n}\Lambda (2n+3,2n+1).
\end{eqnarray*}
We turn to the problem of the evaluations at $(1,-1,0,\ldots ,0).$ In each
case, the value will be expressed in terms of a balanced ${}_{4}F_{3}$%
-series which is obviously positive. This is the result of applying the
Whipple transformation: 
\begin{gather*}
_{4}F_{3}\left( \QATOP{-n,a,b,c}{d,e,f};1\right) = \\
\frac{(1+a-e-n)_{n}(1+a-f-n)_{n}}{(e)_{n}(f)_{n}}{}_{4}F_{3}\left( \QATOP{%
-n,a,d-b,d-c}{d,1+a-e-n,1+a-f-n};1\right)
\end{gather*}
provided $-n+a+b+c+1=d+e+f$ (balanced), and $n\in \mathbf{Z}_{+}$. Setting $%
x=x_{0}=(1,\sqrt{-1},0,\ldots )$ in the basis polynomials produces the
desired values.

\begin{lemma}
For $0\leq j\leq n,$ $\phi _{2n,2j+1}(x_{0})=0,\phi _{2n+1,2j+1}(x_{0})=0$
and 
\begin{eqnarray*}
\phi _{2n,2j}(x_{0}) &=&\frac{(2k+1)_{n+j}(-1)^{n}(-n)_{j}(2k+2n+1)}{n!j!(k+%
\frac{3}{2})_{j}(2k+1)}, \\
\phi _{2n+1,2j}(x_{0}) &=&\frac{2(2k+2)_{n+j}(-1)^{n}(-n)_{j}}{n!j!(k+\frac{3%
}{2})_{j}}.
\end{eqnarray*}
\end{lemma}

\begin{proof}
Substituting $x=x_{0}$ in the generating functions yields $%
\sum_{n=0}^{\infty }\sum_{j=0}^{n}\phi
_{n,j}(x_{0})=(1+2t-t^{2})((1+t^{2})^{2}-4s^{2}t^{2})^{-(k+1)}.$ The latter
term expands to 
\begin{equation*}
\sum_{i,j=0}^{\infty }s^{2j}t^{2i+2j}(-1)^{i}\frac{%
2^{2j}(k+1)_{j}(2k+2+2j)_{i}}{j!i!}.
\end{equation*}
Now multiply top and bottom by $(2k+2)_{2j},$ replace $i+j$ by $n,$ (also a
simple calculation to multiply the resulting series by $(1-t^{2})$), and in
the denominator expand $(2k+2)_{2j}=2^{2j}(k+1)_{j}(k+\frac{3}{2})_{j}.$
\end{proof}

To illustrate the intermediate steps, consider the case 
\begin{eqnarray*}
h_{4n,0}(x_{0}) &=&(-1)^{n}\frac{(2k+1)_{n}(2k+2n+1)}{n!(2k+1)} \\
&&{}_{4}F_{3}\left( \QATOP{-n,k_{2}+(N-1)k+2n,\frac{1}{2},2k+n+1}{k+\frac{3}{%
2},k_{2}+n,Nk+n+1};1\right)
\end{eqnarray*}
and transform the series, using $a=k_{2}+(N-1)k+2n{}$ and $d=k+\frac{3}{2}.$
The other cases are done similarly (the case $h_{4n+3,0}$ incorporates one
additional step, see case \ref{c4n3} above). The results are: 
\begin{eqnarray*}
h_{4n,0}(x_{0}) &=&(-1)^{n}\frac{%
(2k+1)_{n}(2k+2n+1)((N-1)k+1)_{n}(k_{2}-k)_{n}}{%
n!(2k+1)(k_{2}+n)_{n}(Nk+n+1)_{n}} \\
&&{}_{4}F_{3}\left( \QATOP{-n,k_{2}+(N-1)k+2n,k+1,-n-k+\frac{1}{2}}{k+\frac{3%
}{2},k_{2}-k,(N-1)k+1};1\right) , \\
h_{4n+2,0}(x_{0}) &=&(-1)^{n}\frac{2(2k+2)_{n}((N-1)k+1)_{n}(k_{2}-k)_{n}}{%
n!(k_{2}+n+1)_{n}(Nk+n+2)_{n}} \\
&&{}_{4}F_{3}\left( \QATOP{-n,k_{2}+(N-1)k+2n+1,k+1,-n-k-\frac{1}{2}}{k+%
\frac{3}{2},k_{2}-k,(N-1)k+1};1\right) , \\
h_{4n+1,0}(x_{0}) &=&(-1)^{n}\frac{%
(2k+1)_{n}(2k+2n+1)((N-1)k+1)_{n}(k_{2}-k+1)_{n}}{%
n!(2k+1)(k_{2}+n+1)_{n}(Nk+n+1)_{n}} \\
&&{}_{4}F_{3}\left( \QATOP{-n,k_{2}+(N-1)k+2n+1,k+1,-n-k+\frac{1}{2}}{k+%
\frac{3}{2},k_{2}-k+1,(N-1)k+1};1\right) , \\
h_{4n+3,0}(x_{0}) &=&(-1)^{n+1}\frac{%
(2k_{2}+2n+1)(2k+2)_{n}((N-1)k+1)_{n}(k_{2}-k+1)_{n}}{%
n!(k_{2}+n+1)_{n+1}(Nk+n+2)_{n}} \\
&&{}_{4}F_{3}\left( \QATOP{-n,k_{2}+(N-1)k+2n+2,k+1,-n-k-\frac{1}{2}}{k+%
\frac{3}{2},k_{2}-k+1,(N-1)k+1};1\right) .
\end{eqnarray*}
The values of the even parts of $h_{4n+1,1}(x_{0})$ and $h_{4n+3,1}(x_{0})$
are obtained by replacing $k_{2}$ by $k_{2}+1$ in $h_{4n+0,0}(x_{0})$ and $%
h_{4n+2,0}(x_{0})$ respectively (``even part'' refers to $f_{0}$ in the
expressions $h_{2m+1,1}(x)=x_{1}x_{2}f_{0}(x_{1}^{2},x_{2}^{2},\ldots )$).
This completes the calculation of the $L^{2}$-norms of the harmonic
polynomials. \ The $_{4}F_{3}-$series allow no further simplification.

\section{Discussion}

To conclude, we discuss the significance of the results, especially with
regard to applications and indications for further research. The problem
that was solved here is, in a sense, the minimal approach to constructing
harmonic polynomials of type $B$. It may turn out that a different
normalization may be more useful or concise; for example the value $%
h_{4n,0}(1,0,\ldots )$ can be written as $\dfrac{(-1)^{n}(k+1)_{n}\Lambda
(2n,2n)}{n!\Lambda (4n,0)}$ and similar expressions hold for the other
formulae in section \ref{LC}. The expression for $||h_{4n,0}||^{2}$ is also
somewhat simplified by changing the normalization to $\frac{\Lambda (4n,0)}{%
\Lambda (2n,2n)}h_{4n,0}$. Of course the $_{4}F_{3}$ part stays.

\subsection{Application:}

There is a quantum many-body exactly-solvable model associated with $\Delta
_{B}$, namely the spin Calogero model of Yamamoto and Tsuchiya \cite{Y}\cite
{YT}. This deals with $N$ identical particles on a line with inverse-square
mutual repulsion potential and an external harmonic confinement potential.
As well, the particles have a two-valued spin which can be exchanged among
them. The construction of eigenfunctions in terms of nonsymmetric Jack and
generalized Hermite polynomials was discussed in \cite{D6}. The Hamiltonian
for the system ( with $\omega ,k,k_{1}>0)$ is 
\begin{eqnarray*}
\mathcal{H} &=&\sum_{i=1}^{N}\left\{ -\left( \frac{\partial }{\partial x_{i}}%
\right) ^{2}+\omega ^{2}x_{i}^{2}+\frac{k_{1}(k_{1}-\sigma _{i})}{x_{i}^{2}}%
\right\}  \\
&&+2k\sum_{1\leq i<j\leq N}\left\{ \frac{k-\sigma _{ij}}{(x_{i}-x_{j})^{2}}+%
\frac{k-\tau _{ij}}{(x_{i}+x_{j})^{2}}\right\} .
\end{eqnarray*}
The ground state for the system is 
\begin{equation*}
\psi (x)=\prod_{i=1}^{N}|x_{i}|^{k_{1}}\prod_{1\leq i<j\leq
N}|x_{i}^{2}-x_{j}^{2}|^{k}\exp (-\frac{\omega |x|^{2}}{2}).
\end{equation*}

Then the conjugate $\psi \mathcal{H}\psi ^{-1}=2\omega (\sum_{i=1}^{N}x_{i}%
\frac{\partial }{\partial x_{i}}+Nk_{2})-\Delta _{B}.$ Let $f_{m}(x)$ be a
harmonic and homogeneous polynomial of degree $m,$ then for $n=0,1,2,\ldots $
the function $L_{n}^{(c)}(\omega |x|^{2})f_{m}(x)\psi (x)$ is an
eigenfunction of $\mathcal{H}$ with eigenvalue $\allowbreak 2\omega
(m+2n+Nk_{2}),$ where $c=m+Nk_{2}-1.$ Here 
\begin{equation*}
L_{n}^{(c)}(t)=\frac{(c+1)_{n}}{n!}\sum_{j=0}^{n}\frac{(-n)_{i}}{(c+1)_{i}}%
\frac{t^{i}}{i!}
\end{equation*}
denotes the Laguerre polynomial of index $c$ and degree $n.$ The set of all
such functions with $m+2n=s$ spans all the eigenfunctions with eigenvalue $%
2\omega (s+Nk_{2}).$ The set $\{L_{n}^{(c)}(\omega |x|^{2})f_{m}(x)\}$as a
basis for polynomials was used in the study of inner products \cite{D2} and
the Hankel transform \cite{D3}. Van Diejen \cite{vD} considered $W_{N}$%
-invariant eigenfunctions of this type for $\mathcal{H}$. The polynomials $%
h_{2n,0}$ can produce such invariants by summing over translates: 
\begin{equation*}
(1+\sum_{j>2}\sigma _{2j}+\sum_{2<i<j\leq N}\sigma _{1i}\sigma
_{2j})h_{2n,0}(x).
\end{equation*}

\subsection{Further work:}

It is still an open problem to find an orthogonal basis for the harmonic
homogeneous polynomials. Such bases are useful in approximation theory and
numerical cubature (see Xu \cite{X1}\cite{X2}). It is not difficult to write
down self-adjoint operators on polynomials, for example $%
(x_{i}T_{j}-x_{j}T_{i})^{2}$ for $1\leq i<j\leq N.$ This is a useful method
for abelian reflection groups. However computer algebra calculations reveal
that the characteristic polynomials of these operators on polynomials of not
large degree do not factor linearly in $\mathbf{Q}(k,k_{1})$ (for type $B$).
Hence one does not expect tractable eigenfunction decompositions. It seems
worthwhile to try to extend to more variables the generating function
construction for $F_{0},F_{1}$ which was the main device for this paper
(that is, consider harmonic polynomials annihilated by $T_{i}$ for $i>n_{0}$%
; already the case $n_{0}=3$ is interesting). Obviously there will need to
be a more sophisticated way of handling the different cases. The present
approach is just barely tolerable for the different basis functions involved
in the representation theory of $B_{2}.$ It certainly seems that finding
orthogonal bases is considerably more complicated then the construction of
nonsymmetric Jack polynomials.

\end{document}